# On probabilities associated with the minimum distance between events of a Poisson process in a finite interval

Shai Covo
Department of Mathematics
Bar Ilan University
52900 Ramat-Gan, Israel
E-mail: green355@netvision.net.il

November 1, 2018

**Abstract**

We revisit the probability that any two consecutive events in a Poisson process $N$ on $[0,t]$ are separated by a time interval which is greater than $s\,(<t)$ (a particular scan statistic probability), and the closely related probability (recently introduced by Todinov [8], who denotes it as $p_{\mathrm{MFFOP}}$) that before any event of $N$ in $[0,t]$ there exists an event-free interval greater than $s$. Both probabilities admit simple explicit expressions, which, however, become intractable for very large values of $t/s$. Our main objective is to demonstrate that these probabilities can be approximated extremely well for large values of $t/s$ by some very tractable and attractive expressions (actually, already for $t$ larger than a few multiples of $s$).

## 1 Introduction

Throughout this paper, $N = \{N_t : t \geq 0\}$ denotes a Poisson process with rate $\lambda > 0$. For any $s > 0$ and $t > 0$, we define

$$\phi(s,t;\lambda) = \begin{cases} \mathrm{P}(N_t < 2) \ (= \mathrm{e}^{-\lambda t}(1+\lambda t)) & , \ 0 < t \leq s; \\ \mathrm{P}(\max_{0 \leq r \leq t-s}[N_{r+s} - N_r] < 2) & , \ t > s. \end{cases}$$

By convention, $\phi(s, t; \lambda) = 1$ if $t = 0$. That is, $\phi(s, t; \lambda)$, for $s, t > 0$, is the probability that $N$ has no two consecutive events in $[0, t]$ separated by a time interval which is shorter than $s$. Its complement is denoted by

$$\bar{\phi}(s, t; \lambda) := 1 - \phi(s, t; \lambda),$$

and gives us the distribution function of the random variable $\tau(s; \lambda)$ defined to be the waiting time until two events of $N$ occur within $s$ time units; namely,

$$\mathrm{P}(\tau(s; \lambda) \leq t) = \bar{\phi}(s, t; \lambda), \quad t \geq 0. \tag{1.1}$$

As indicated in the abstract, the probability $\phi(s, t; \lambda)$ is a particular *scan statistic* probability. Indeed, the random variable $Y_N(s, t)$ defined by $Y_N(s, t) = \max_{0 \leq r \leq t-s}[N_{r+s} - N_r]$ is the scan statistic of 'window size' $s$ for the Poisson process $N$ on $[0, t]$. The distribution of $Y_N(s, t)$ has been extensively studied in the literature. However, the particular scan statistic probability $\mathrm{P}(Y_N(s, t) < 2)$ seems to have received little attention in the literature relative to its importance (and simplicity). Yet, as the author discovered after most of this research was completed, the fundamental result of this paper has been essentially identified in the (scan statistic) literature (Gates and Westcott [3]); nevertheless, our overall presentation of the result has several advantages and novelties. This point will be discussed in some detail in Appendix A, thus highlighting the contribution of our paper. It should also be stressed here that approximations to the general scan statistic probability $\mathrm{P}(Y_N(s, t) < k)$, $k \geq 2$, which are readily available in the literature, are usually not satisfactory when applied to our case of interest, namely, $k = 2$. A familiar and relatively successful one is due to Naus [5]; however, it is far too inferior compared to our results—both from the accuracy point of view, and from the simplicity one (see also Appendix A).

Closely related to the probability $\phi(s, t; \lambda)$ is the probability that before any event of $N$ in $[0, t]$ there exists an event-free interval greater than $s$. We shall denote it by $\phi_2(s, t; \lambda)$. This probability has recently been introduced and studied by Todinov [8], [9], [10] (to cite a few) in the context of reliability engineering (sometimes in a more general setting), where it is denoted as $p_{\mathrm{MFFOP}}$. The acronym MFFOP stands for "Minimum failure-free operating period"; thus, in Todinov's framework, the events of $N$ are interpreted as "failures".

It follows straight from the definitions that

$$\phi_2(s, t; \lambda) = \mathrm{e}^{-\lambda s} \phi(s, t - s; \lambda), \quad t > s. \tag{1.2}$$

Indeed, the event associated with the probability $\phi_2(s, t; \lambda)$ can occur if and only if $N_s = 0$ and $N$ has no two consecutive events in $[s, t]$ separated by a time interval which is shorter than $s$. Relation (1.2) allows us to focus on the probability $\phi(s, t; \lambda)$ (which is more interesting, in our view).

The probability $\phi(s, t; \lambda)$, for $s, t > 0$, admits the following expression:

$$\phi(s, t; \lambda) = e^{-\lambda t}\left(1 + \sum_{k=1}^{\lceil t/s \rceil} \frac{\lambda^k [t - (k-1)s]^k}{k!}\right) \quad (1.3)$$

(Todinov [8, Eq. (17)], [9, Eq. (4)], [10, Eq. (15.4)]), where $\lceil \cdot \rceil$ denotes the ceiling function. For completeness of exposition and convenience of the reader, we give here two full proofs of this simple result.

Writing the right-hand side of (1.3) as

$$e^{-\lambda t}(1 + \lambda t) + \sum_{k=2}^{\infty} \frac{e^{-\lambda t}(\lambda t)^k}{k!}\left(1 - \frac{(k-1)s}{t}\right)^k \mathbf{1}\left(\frac{t}{k-1} > s\right),$$

where $\mathbf{1}(\cdot)$ denotes the indicator function, and noting that the greatest possible minimum distance between $k \geq 2$ points in $[0, t]$ is $t/(k-1)$, it follows by the law of total probability (conditioning on the number of events of $N$ in $[0, t]$) that (1.3) is true if $(1 - (k-1)s/t)^k$ is the probability that the minimum distance between $k \geq 2$ i.i.d. uniform random variables on $[0, t]$ is greater than $s$. The latter is a well-known fact; however, we take this opportunity to prove the following more general fact (Todinov [8, Eq. (7)]).

Let $U_{(1)}, \ldots, U_{(n)}$ be $n$ order statistics from a uniform distribution on $[0, t]$, and $s_{1,2}, \ldots, s_{n-1,n}$ nonnegative numbers with $\sum_{i=1}^{n-1} s_{i,i+1} \leq t$. Then

$$P(U_{(i+1)} - U_{(i)} > s_{i,i+1} \; \forall i = 1, \ldots, n-1) = \left(1 - \frac{1}{t}\sum_{i=1}^{n-1} s_{i,i+1}\right)^n. \quad (1.4)$$

The proof we now give is a rigorous justification of Todinov's explanation in [7, Section 3.1], and a short, straightforward alternative to his derivation in [8]. Define

$$A = \{(u_1, \ldots, u_n) \in [0, t]^n : u_{i+1} - u_i > s_{i,i+1} \; \forall i = 1, \ldots, n-1\}$$

and

$$B = \{(u_1, u_2 - s_{1,2}, \ldots, u_n - (s_{1,2} + \cdots + s_{n-1,n})) : (u_1, \ldots, u_n) \in A\} \\ = \{(v_1, v_2, \ldots, v_n) : 0 \leq v_1 < \cdots < v_n \leq t - (s_{1,2} + \cdots + s_{n-1,n})\}. \quad (1.5)$$

Recall that the joint distribution of the order statistics in (1.4) is given by $P(U_{(1)} \in du_1, \ldots, U_{(n)} \in du_n) = n!/t^n \, du_1 \cdots du_n$, $0 \leq u_1 < \cdots < u_n \leq t$. Thus, the left-hand side of (1.4) is equal to $n! t^{-n} \int_A du_1 \cdots du_n$, and hence also to $n!t^{-n} \int_B du_1 \cdots du_n$ since, as is clear from the first equality in (1.5), $A$ and $B$ have the same Lebesgue measure. On the other hand, by the second equality in (1.5), $n!t^{-n} \int_B du_1 \cdots du_n$ is the probability of $\{U_{(n)} \leq t - \sum_{i=1}^{n-1} s_{i,i+1}\}$. Thus, (1.4) is established.

We now give an alternative, original proof of (1.3). The 'base case' $\lceil t/s \rceil = 1$ (i.e., $0 < t \leq s$) holds by definition. Assume that $t > s$. By conditioning on the time of the first event of $N$, we obtain

$$\phi(s, t; \lambda) = \int_0^{t-s} e^{-\lambda s} \phi(s, t - (u+s); \lambda) \lambda e^{-\lambda u} \, du$$
$$+ \int_{t-s}^t e^{-\lambda(t-u)} \lambda e^{-\lambda u} \, du$$
$$+ \int_t^\infty \lambda e^{-\lambda u} \, du,$$

leading to

$$\phi(s, t; \lambda) = \begin{cases} e^{-\lambda t}(1 + \lambda t) & , 0 < t \leq s; \\ e^{-\lambda t}(1 + \lambda s) + \int_s^t \lambda e^{-\lambda u} \phi(s, t-u; \lambda) \, du & , t > s. \end{cases} \quad (1.6)$$

The rest of the proof is by induction on $\lceil t/s \rceil = n$. Suppose that (1.3) holds for $n$, and let $\lceil t/s \rceil = n+1$. Making the change of variable $t - u \mapsto u$ in (1.6) and then applying the induction hypothesis, we obtain

$$\begin{aligned}
\phi(s, t; \lambda) &= e^{-\lambda t}(1 + \lambda s) + e^{-\lambda t} \int_0^{t-s} \lambda e^{\lambda u} \phi(s, u; \lambda) \, du \\
&= e^{-\lambda t}(1 + \lambda s) + e^{-\lambda t} \int_0^{t-s} \lambda \left(1 + \sum_{j=1}^{\lceil u/s \rceil} \frac{\lambda^j [u - (j-1)s]^j}{j!}\right) du \\
&= e^{-\lambda t}(1 + \lambda t) + e^{-\lambda t} \lambda \sum_{j=1}^n \int_{(j-1)s}^{t-s} \frac{\lambda^j [u - (j-1)s]^j}{j!} \, du \\
&= e^{-\lambda t} \left(1 + \sum_{k=1}^{n+1} \frac{\lambda^k [t - (k-1)s]^k}{k!}\right).
\end{aligned} \quad (1.7)$$

Thus, (1.3) is established.

*Remark* 1.1. It is worth noting that $\phi$ depends on $s$, $t$, and $\lambda$ only through $t/s$ and $\lambda s$. In particular, we have that $\phi(s, t; \lambda) = \phi(1, t/s; \lambda s)$.

Taking into account the convolution theorem for Laplace transforms, Eq. (1.6) invites us to calculate the Laplace transform of $\phi(s, t; \lambda)$ as a function of $t$ (for $s$ fixed). Anticipating Theorem 2.2 below, we note that, based on (1.3), it can be shown that the equality

$$\int_0^\infty e^{-ut}\phi(s,t;\lambda)\,dt = \frac{1}{u+\lambda} + \frac{\lambda}{u+\lambda}\cdot\frac{1}{u+\lambda-\lambda e^{-(u+\lambda)s}} \qquad (1.8)$$

holds for all $u > u_c := -\lambda[1 - B(\lambda s)]$, where $B(\lambda s)$ is the solution $B \in (0,1)$ of $e^{-\lambda sB} = B$, while for $u \leq u_c$ the integral diverges. This number $u_c$ is thus the *abscissa of convergence* of the Laplace transform of $\phi(s,t;\lambda)$; it is, trivially, also the abscissa of convergence of the Laplace transform of $Ae^{u_ct}$, where $A = A(\lambda, s) > 0$ is a constant with respect to $t$. Theorem 2.2 gives the precise relation between the functions $\phi(s,t;\lambda)$ and $e^{u_ct}$. Actually, the origin of Theorem 2.2 lies in the following remark.

*Remark* 1.2. It follows from either (1.3) or the first equality in (1.7) that

$$\frac{\partial}{\partial t}\phi(s,t;\lambda) = -\lambda\phi(s,t;\lambda) + \lambda e^{-\lambda s}\phi(s,t-s;\lambda), \quad t > s.$$

With the above notation, we note that this differential equation is also satisfied by the function $Ae^{u_ct}$.

Recall (1.2). Then, from (1.3) it follows immediately that, for $s, t > 0$,

$$\phi_2(s,t;\lambda) = e^{-\lambda t}\left(1 + \sum_{k=1}^{\lceil t/s \rceil - 1} \frac{\lambda^k(t-ks)^k}{k!}\right), \qquad (1.9)$$

which is Todinov's formula, given e.g. in [8, Eq. (16)], [9, Eq. (3)], and [10, Eq. (15.2)]. (However, relation (1.2) has not been observed in Todinov's works.) Just as in Remark 1.1, we have that $\phi_2(s,t;\lambda) = \phi_2(1, t/s; \lambda s)$.

To conclude the introduction, we note the following simple facts.

Let $\xi_1, \xi_2, \ldots$ be the interarrival times of $N$, and set $Y = \inf\{i : \xi_i < s\}$. Then, $Y$ is a geometric random variable with $P(Y = i) = e^{-\lambda s(i-1)}(1 - e^{-\lambda s})$, and the following equality in distribution holds:

$$\tau_2(s;\lambda) := \sum_{i=1}^Y \xi_i \stackrel{d}{=} \sum_{i=1}^{Y-1}(s + \tilde{\xi}_i) + \zeta, \qquad (1.10)$$

where $\tilde{\xi}_i$ (just as $\xi_i$) are independent exponentials with mean $1/\lambda$, independent also of $Y$ and of the random variable $\zeta$ with density $(1 - e^{-\lambda s})^{-1}\lambda e^{-\lambda x}$, $0 < x < s$. The random variable $\tau_2$ is related to the probability $\phi_2$ by

$$P(\tau_2(s;\lambda) \leq t) = 1 - \phi_2(s,t;\lambda), \quad t \geq 0.$$

Using the decomposition on the right-hand side of (1.10), one can show that

$$E[\tau_2(s;\lambda)] = \frac{1}{\lambda(1 - e^{-\lambda s})} \quad \text{and} \quad \text{Var}[\tau_2(s;\lambda)] = \frac{1 + 2\lambda s e^{-\lambda s}}{\lambda^2(1 - e^{-\lambda s})^2}.$$

(Actually, the formula for the expectation follows immediately from Wald's identity, applied to $\tau_2(s;\lambda) = \sum_{i=1}^{Y} \xi_i$ .) Then, from the obvious relation

$$\tau(s;\lambda) \stackrel{\mathrm{d}}{=} \xi + \tau_2(s;\lambda),$$

where $\xi$ is exponential with mean $1/\lambda$, independent of $\tau_2(s;\lambda)$, we immediately obtain the expectation and variance of $\tau(s;\lambda)$.

## 2  Main results

Although (1.3) gives us a neat result, it may not be effective for very large values of $t/s$. This point will be discussed in the next section. Our first main result, Theorem 2.1 below, gives an asymptotic power series expansion for $\phi(s,t;\lambda)$ as $s \downarrow 0$. Before we state and prove it, we introduce the following polynomials, defined for $x > 0$:

$$P_n(x) = \sum_{k=0}^{\infty} \frac{\mathrm{e}^{-x} x^k}{k!} [k + (n-1)]^n$$

$$= \sum_{i=0}^{n} \binom{n}{i} (n-1)^{n-i} \sum_{j=0}^{i} S(i,j) x^j,$$

where $S(i,j)$ denote the Stirling numbers of the second kind and where we have used the well-known identity

$$\sum_{k=0}^{\infty} \frac{\mathrm{e}^{-x} x^k}{k!} k^i = \sum_{j=0}^{i} S(i,j) x^j, \quad i \in \mathbb{Z}_+,$$

expressing the moments of the Poisson distribution with mean $x$ in terms of Stirling numbers of the second kind (see e.g. Comtet [2, p. 160]). The $S(i,j)$ can be calculated recursively as follows (see e.g. Comtet [2, p. 208]):

$$S(i,j) = S(i-1, j-1) + jS(i-1, j) \quad \text{for } 2 \le j < i,$$

with $S(i,1) = S(i,i) = 1$ for $i \ge 1$. In addition, $S(i,0) := \mathbf{1}_{\{0\}}(i)$. The first eight polynomials $P_n(x)$ are listed in Table 1.

**Theorem 2.1.** *Fix $t > 0$. Then,*

$$\phi(s,t;\lambda) \sim \sum_{n=0}^{\infty} (-1)^n P_n(\lambda t) \frac{\lambda^n s^n}{n!} \quad \text{as } s \downarrow 0.$$

*That is, for each $N = 0, 1, 2, \ldots$,*

$$\phi(s,t;\lambda) - \sum_{n=0}^{N} (-1)^n P_n(\lambda t) \frac{\lambda^n s^n}{n!} = O(s^{N+1}) \quad \text{as } s \downarrow 0. \tag{2.1}$$

Table 1: The first eight polynomials $P_n(x)$.

| $n$ | $P_n(x)$ |
|---|---|
| 0 | 1 |
| 1 | $x$ |
| 2 | $1 + 3x + x^2$ |
| 3 | $8 + 19x + 9x^2 + x^3$ |
| 4 | $81 + 175x + 97x^2 + 18x^3 + x^4$ |
| 5 | $1024 + 2101x + 1275x^2 + 305x^3 + 30x^4 + x^5$ |
| 6 | $15625 + 31031x + 19981x^2 + 5590x^3 + 740x^4 + 45x^5 + x^6$ |
| 7 | $279936 + 543607x + 365001x^2 + 113701x^3 + 18200x^4 + 1526x^5 + 63x^6 + x^7$ |

*Proof.* We will show that (2.1) holds for all $N = 1, 3, 5, \ldots$, the case $N = 0, 2, 4, \ldots$ being analogous. Fix $t > 0$ and $N \in \{1, 3, 5, \ldots\}$, and let $K := \lceil t/s \rceil$ be arbitrarily large. It follows from (1.3), using a simple property of the binomial expansion, that $\phi(s, t; \lambda)$ is bounded from above by

$$e^{-\lambda t}\left\{1 + \lambda t + \sum_{k=2}^{N+1} \frac{\lambda^k [t - (k-1)s]^k}{k!} + \sum_{k=N+2}^{K} \frac{\lambda^k}{k!} \sum_{i=0}^{N+1} (-1)^i \binom{k}{i} t^{k-i} [(k-1)s]^i \right\}. \tag{2.2}$$

Collecting powers of $s$, it follows straightforwardly that (2.2) is equal to

$$\sum_{n=0}^{N+1} (-1)^n \frac{\lambda^n s^n}{n!} \left\{ P_n(\lambda t) - \sum_{k=K-n+1}^{\infty} \frac{e^{-\lambda t}(\lambda t)^k}{k!} [k + (n-1)]^n \right\}. \tag{2.3}$$

It is easy to show that the second term in the curly brackets is $o(s^r)$ as $s \downarrow 0$, for any $r > 0$. We conclude that the left-hand side of (2.1) is bounded from above by $|O(s^{N+1})|$ as $s \downarrow 0$. The analogue of (2.2) for a lower bound is obtained by replacing $N+1$ by $N$ and $N+2$ by $N+1$; continuing as in the upper-bound case, we conclude that the left-hand side of (2.1) is bounded from below by $o(s^r)$ as $s \downarrow 0$, for any $r > 0$. Thus (2.1) is proved (for $N$ odd, and analogously for $N$ even). ∎

*Remark* 2.1. From the proof of Theorem 2.1, we can obtain the following result:

$$\phi(s, t; \lambda) \to e^{-\lambda^2 c} \quad \text{as } s \downarrow 0, \ t \to \infty, \ st \to c, \tag{2.4}$$

for $c > 0$ fixed. Indeed, with $s$ and $t$ as in (2.4), $P_n(\lambda t)\lambda^n s^n \to (\lambda^2 c)^n$ (for all fixed $n$) and, as is easy to show, the second term in the curly brackets of (2.3) tends to 0. It is straightforward to verify that this result is consistent with the approximation (2.16), below.

*Remark* 2.2. Similarly to Remark 2.1, from the proof of Theorem 2.1, we can obtain the following result:

$$\phi(s,t;\lambda) \to e^{-ct} \quad \text{as } s \downarrow 0, \ \lambda \to \infty, \ \lambda^2 s \to c, \tag{2.5}$$

for $c > 0$ fixed. Indeed, with $s$ and $\lambda$ as in (2.5), $P_n(\lambda t)\lambda^n s^n \to (ct)^n$ (for all fixed $n$) and, as is easy to show, the second term in the curly brackets of (2.3) tends to 0. (Actually, (2.5) can be deduced from (2.4) by applying the relation $\phi(s,t;\lambda) = \phi(1,t/s;\lambda s)$.) Again, it is straightforward to verify that this result is consistent with the approximation (2.16), below. Now, recalling (1.1), the following statement is immediate. With $s$ and $\lambda$ as in (2.5), the waiting time until two events of $N$ occur within $s$ time units converges in distribution to an exponential random variable with mean $1/c$.

Our second main result, Theorem 2.2 below, provides a tractable asymptotic expression for $\phi(s,t;\lambda)$ as $t \to \infty$. As we shall confirm numerically in the next section, this expression turns out to be an extremely successful approximation to $\phi(s,t;\lambda)$ for large values of $t/s$ in general—far more successful than one would expect a priori; surprisingly, it is even superior to its counterpart (2.1) for the case when $t$ is fixed and $s$ is small.

We first introduce some more notation. We let $\varphi(s,t;\lambda)$ denote the function defined by $\varphi(s,t;\lambda) = e^{\lambda t}\phi(s,t;\lambda)$, $s > 0$, $t > 0$. By Remark 1.1, the function $\varphi(s,t;\lambda)$ is a function of $x := \lambda s$ and $r := t/s$; we denote this function of $x$ and $r$ by $\varphi(x,r)$. Specifically, for any $x > 0$ and $r > 0$,

$$\varphi(x,r) = 1 + \sum_{k=1}^{[r]} \frac{x^k[r-(k-1)]^k}{k!}. \tag{2.6}$$

We further denote by $B(x)$ the function $B : (0,\infty) \to (0,1)$ defined as the solution of the implicit equation

$$\exp[-xB(x)] = B(x),$$

and by $A(x)$ the function defined by

$$A(x) = \frac{1}{B(x)[1+xB(x)]}. \tag{2.7}$$

The functions $A(x)$ and $B(x)$ are plotted in Figure 1. Finally, we let $A(x,r)$ denote the function defined by

$$A(x,r) = \frac{\varphi(x,r)}{\exp[xB(x)r]}. \tag{2.8}$$

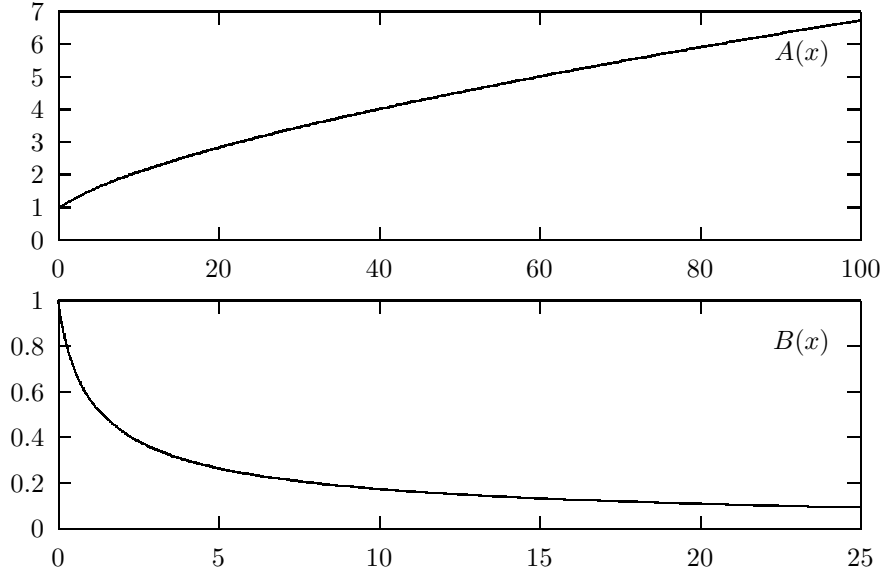

Figure 1: The functions $A(x)$ and $B(x)$.

Thus, $\phi(s, t; \lambda)$ can be written as

$$\phi(s, t; \lambda) = A(\lambda s, t/s)e^{-\lambda[1 - B(\lambda s)]t}. \tag{2.9}$$

As is customary, $f \sim g$ means that $f/g \to 1$.

**Theorem 2.2.** *The following convergence holds for any $x > 0$:*

$$A(x, r) \to A(x) \quad \text{as } r \to \infty. \tag{2.10}$$

*Hence, for any fixed $s > 0$,*

$$\phi(s, t; \lambda) \sim A(\lambda s)e^{-\lambda[1 - B(\lambda s)]t} \quad \text{as } t \to \infty. \tag{2.11}$$

*Proof.* To prove the theorem, we will use a well-known result from analysis, the Wiener–Ikehara theorem (see Appendix B). In order to apply it, we first rewrite (2.10) as

$$e^{-r}\varphi(x, \tfrac{r}{xB(x)}) \to A(x) \quad \text{as } r \to \infty. \tag{2.12}$$

Next, we have to calculate the Laplace–Stieltjes transform of $\varphi(x, \tfrac{r}{xB(x)})$ as a function of $r$ (for $x$ fixed). We denote it by $f(z; x)$, where $z = u + iv$ with $u > 1$. We claim that

$$f(z; x) = 1 + \frac{1}{zB(x) - \exp[-zxB(x)]}, \quad u > 1. \tag{2.13}$$

As a first step, let us justify the following equality:

$$f(z;x) := \int_{0^-}^{\infty} e^{-zr}\,d\varphi(x, \tfrac{r}{xB(x)}) = z\int_0^{\infty} e^{-zr}\varphi(x, \tfrac{r}{xB(x)})\,dr,$$

under the assumption that the right-hand side exists if $u > 1$. Integrating by parts, we have

$$\int_{0^-}^{M} e^{-zr}\,d\varphi(x, \tfrac{r}{xB(x)}) = e^{-zM}\varphi(x, \tfrac{M}{xB(x)}) + z\int_0^{M} e^{-zr}\varphi(x, \tfrac{r}{xB(x)})\,dr.$$

We thus have to show that $e^{-zM}\varphi(x, \tfrac{M}{xB(x)}) \to 0$ as $M \to \infty$. This is easily seen to hold if $z$ is real ($> 1$), for the left-hand side is then increasing in $M$ and hence converges. Hence, it also holds for complex $z$, by taking absolute value. It now follows that (2.13) will be proved if we show that

$$z\lim_{\substack{n\to\infty\\n\in\mathbb{N}}} \int_0^{n} e^{-zr}\varphi(x, \tfrac{r}{xB(x)})\,dr = 1 + \frac{1}{zB(x) - \exp[-zxB(x)]}, \quad u > 1 \quad (2.14)$$

(consider first real $z$, then complex). Making the change of variable $r/(xB(x)) \mapsto r$ and then substituting from (2.6), we obtain

$$\int_0^{n} e^{-zr}\varphi(x, \tfrac{r}{xB(x)})\,dr$$

$$= xB(x)\sum_{j=1}^{n}\int_{j-1}^{j} e^{-zxB(x)r}\left(1 + \sum_{k=1}^{j} \frac{x^k[r-(k-1)]^k}{k!}\right)dr$$

$$= xB(x)\left\{\int_0^{n} e^{-zxB(x)r}\,dr + \sum_{j=1}^{n}\int_{j-1}^{n} e^{-zxB(x)r}\frac{x^j[r-(j-1)]^j}{j!}\,dr\right\}$$

$$= xB(x)\left\{\int_0^{n} e^{-zxB(x)r}\,dr + \sum_{j=1}^{n}\int_0^{n+1-j} e^{-zxB(x)(r+j-1)}\frac{x^j r^j}{j!}\,dr\right\}.$$

We now fix $M \in \mathbb{N}$ (arbitrarily large) and let $n > M$, to be sent to $\infty$. A straightforward calculation, using the gamma function, gives

$$\left|\sum_{j=M}^{n}\int_0^{n+1-j} e^{-zxB(x)(r+j-1)}\frac{x^j r^j}{j!}\,dr\right| \leq \frac{1}{x[uB(x)]^2}\sum_{j=M-1}^{n-1}\left[\frac{e^{-uxB(x)}}{uB(x)}\right]^j.$$

Moreover, it follows immediately from the definition of $B(x)$ that $e^{-uxB(x)} < uB(x)$ (for $u > 1$). Hence, from the arbitrariness of $M$, we conclude that

$$\lim_{n\to\infty}\int_0^{n} e^{-zr}\varphi(x, \tfrac{r}{xB(x)})\,dr$$

$$= xB(x)\left\{\int_0^{\infty} e^{-zxB(x)r}\,dr + \sum_{j=1}^{\infty}\int_0^{\infty} e^{-zxB(x)(r+j-1)}\frac{x^j r^j}{j!}\,dr\right\},$$

provided that the infinite sum converges; indeed, it converges whenever $|e^{-zxB(x)}/(zB(x))| < 1$, which is clearly implied by $u > 1$. The last assertion follows from the equality

$$\int_0^\infty e^{-zxB(x)(r+j-1)} \frac{x^j r^j}{j!} \, dr = \frac{1}{x[zB(x)]^2} \left[\frac{e^{-zxB(x)}}{zB(x)}\right]^{j-1},$$

which, in turn, leads straight to (2.14). Hence (2.13) is proved. Thus, in order to prove (2.12) (and hence the theorem) using the Wiener–Ikehara theorem, it suffices to show that the analytic function

$$g(z;x) = 1 + \frac{1}{zB(x) - \exp[-zxB(x)]} - \frac{A(x)}{z-1}, \quad u > 1, \qquad (2.15)$$

has a boundary function $g(1 + iv; x)$ in the following sense. For $u \downarrow 1$, the function $g_u(iv; x) := g(u + iv; x)$ converges to $g(1 + iv; x)$ uniformly on every finite interval $\{-l < v < l\}$. This is clearly satisfied for $\{\delta < |v| < l\}$, for arbitrarily small $\delta > 0$. The uniform convergence on $\{-l < v < l\}$ will thus be established upon showing that $g(u + iv; x)$ has a limit as $u \downarrow 1$ and $v \to 0$. In fact, reducing the fractions in (2.15) to a common denominator and applying l'Hôpital's rule twice (at $z = 1$), we find that

$$\lim_{u \downarrow 1,\, v \to 0} g(u + iv; x) = 1 + \frac{x^2 B(x)}{2[1 + xB(x)]^2}.$$

The theorem is thus established. ∎

It is instructive to note here the following remark.

*Remark* 2.3. Suppose that $\hat{h}(u) := u \int_0^\infty e^{-ur} h(r) \, dr$—the 'Laplace–Stieltjes transform' of a function $h$—exists for $u > 0$, and that $\lim_{u \downarrow 0} \hat{h}(u) = c$. Under a suitable condition (*Tauberian condition*) on $h$, $\lim_{r \to \infty} h(r) = c$ (Final Value Theorem). When applied to the function $A(x, r)$, $x$ fixed, this yields (2.10) (as follows readily from (1.8)). However, it might be difficult to show that $A(x, r)$ fulfills a required Tauberian condition for the validity of this application; cf., e.g., Theorem 1.7.6 in Bingham, Goldie, & Teugels [1] (Karamata's Tauberian Theorem; extended form). The Wiener–Ikehara theorem, which is a complex Tauberian theorem, provides a convenient way to establish the desired result; an alternative is indicated in the last paragraph of Appendix A.

Since the function $\varphi(x, r)$ plays a crucial role in this paper, it is worth making the following easy observation.

*Remark* 2.4. The function $\varphi(x, r)$—as well as its asymptotic expression $A(x)e^{xB(x)r}$ as $r \to \infty$—satisfies the differential equation

$$\frac{\partial}{\partial r} \varphi(x, r) = x\varphi(x, r-1), \quad r > 1.$$

In view of (2.9) and Theorem 2.2, it is natural to define and study the following approximation, which is designed for large values of $t/s$:

$$\phi^\star(s,t;\lambda) := A(\lambda s)e^{-\lambda[1-B(\lambda s)]t} \approx \phi(s,t;\lambda). \quad (2.16)$$

(It indeed deserves to be denoted with a '$\star$'.) As we shall see in the next section, the approximation (2.16) may be effective already for $t$ larger than a few multiples of $s$, and, moreover, even in the entire range of interest $\{t \geq s\}$.

To conclude this section, recall the function $\phi_2(s,t;\lambda)$ given by (1.9). From (2.16) and relation (1.2), it follows straightforwardly the following approximation (appropriate for $s$ and $t$ as above):

$$\phi_2(s,t;\lambda) \approx A(\lambda s)B(\lambda s)e^{-\lambda[1-B(\lambda s)]t}$$
$$= \frac{1}{1+\lambda s B(\lambda s)}e^{-\lambda[1-B(\lambda s)]t}.$$

This establishes the interesting relation $\phi_2(s,t;\lambda) \approx B(\lambda s)\phi(s,t;\lambda)$.

## 3 Numerical results

In this section we demonstrate, in various ways, how extremely successful the approximation (2.16) is. We first note that since

$$\phi^\star(s,t;\lambda) - \phi(s,t;\lambda) = [A(\lambda s) - A(\lambda s, t/s)]e^{-\lambda[1-B(\lambda s)]t}, \quad (3.1)$$

we will pay particular attention to the difference between $A(x)$ and $A(x,r)$. However, the following point should be kept in mind in the sequel. On the one hand, as is seen from Table 2 and Figure 2, the quality of the approximation $A(x,r) \approx A(x)$ reduces as $x$ increases. On the other hand, the larger $x\,(=\lambda s)$ is, the closer the negative exponential term in (3.1) will be to 0.

We begin by considering the following example, referring to (2.4).

*Example* 3.1. We have implemented the computation of $\phi(s,t;\lambda)$ (Eq. (1.3)) and $\phi^\star(s,t;\lambda)$ (Eq. (2.16)) in the Java language. Since, as we shall see later, $\phi(s,t;\lambda)$ and $\phi^\star(s,t;\lambda)$ are extremely close to each other for large values of $t/s$, one may study the convergence (2.4) based on $\phi^\star$ values. As a concrete example, we let $\lambda = 1$, $s = 1/\sqrt{r}$, and $t = \sqrt{r}$, for $r$ very large. By (2.4), $\phi(s,t;\lambda) \to e^{-1}$ as $r \to \infty$. When $r$ was set to $m \cdot 10^4$, for $m = 8, 9$, and $10$, the computational time for $\phi$ was too long: approximately $11.1, 14.0$, and $17.3$ minutes, respectively. The $\phi^\star$ values, on the other hand, were computed in 'no time'. Computation of $\phi^\star(s,t;\lambda) - e^{-1}$ for various values of $r$ indicated that this difference is asymptotically $1.5\mathrm{e}^{-1}/\sqrt{r}$ as $r \to \infty$; this was confirmed by more accurate calculations, suitable for $r$ as large as $10^{60}$, using WIMS function calculator (see below).

Table 2: Tabulated rounded values of $A(x,r) - A(x)$ for various values of $x$ and $r$.

| $x$ | $r$ | | | |
|---|---|---|---|---|
| | 10 | 25 | 50 | 100 |
| 0.010 | $+6 \cdot 10^{-33}$ | $-7 \cdot 10^{-75}$ | $+1 \cdot 10^{-148}$ | $+9 \cdot 10^{-293}$ |
| 0.025 | $+3 \cdot 10^{-28}$ | $+1 \cdot 10^{-64}$ | $+9 \cdot 10^{-126}$ | $-1 \cdot 10^{-247}$ |
| 0.050 | $-7 \cdot 10^{-25}$ | $-4 \cdot 10^{-56}$ | $-9 \cdot 10^{-109}$ | $+3 \cdot 10^{-214}$ |
| 0.100 | $-9 \cdot 10^{-21}$ | $+7 \cdot 10^{-48}$ | $-2 \cdot 10^{-92}$ | $-9 \cdot 10^{-182}$ |
| 0.250 | $+9 \cdot 10^{-17}$ | $-3 \cdot 10^{-37}$ | $+4 \cdot 10^{-73}$ | $+2 \cdot 10^{-142}$ |
| 0.500 | $+2 \cdot 10^{-13}$ | $+1 \cdot 10^{-30}$ | $+9 \cdot 10^{-59}$ | $-2 \cdot 10^{-115}$ |
| 1.000 | $-6 \cdot 10^{-11}$ | $+8 \cdot 10^{-25}$ | $+9 \cdot 10^{-49}$ | $+5 \cdot 10^{-93}$ |
| 2.500 | $+7 \cdot 10^{-9}$ | $-5 \cdot 10^{-19}$ | $+2 \cdot 10^{-35}$ | $-8 \cdot 10^{-70}$ |
| 5.000 | $+1 \cdot 10^{-6}$ | $-4 \cdot 10^{-15}$ | $+8 \cdot 10^{-29}$ | $-4 \cdot 10^{-57}$ |
| 10.00 | $-9 \cdot 10^{-6}$ | $+4 \cdot 10^{-12}$ | $+3 \cdot 10^{-24}$ | $+8 \cdot 10^{-46}$ |
| 25.00 | $-5 \cdot 10^{-4}$ | $-3 \cdot 10^{-9}$ | $-7 \cdot 10^{-18}$ | $-3 \cdot 10^{-35}$ |
| 50.00 | $+1 \cdot 10^{-3}$ | $+2 \cdot 10^{-7}$ | $-7 \cdot 10^{-15}$ | $+3 \cdot 10^{-29}$ |
| 100.0 | $+2 \cdot 10^{-2}$ | $-4 \cdot 10^{-6}$ | $+2 \cdot 10^{-12}$ | $-2 \cdot 10^{-24}$ |

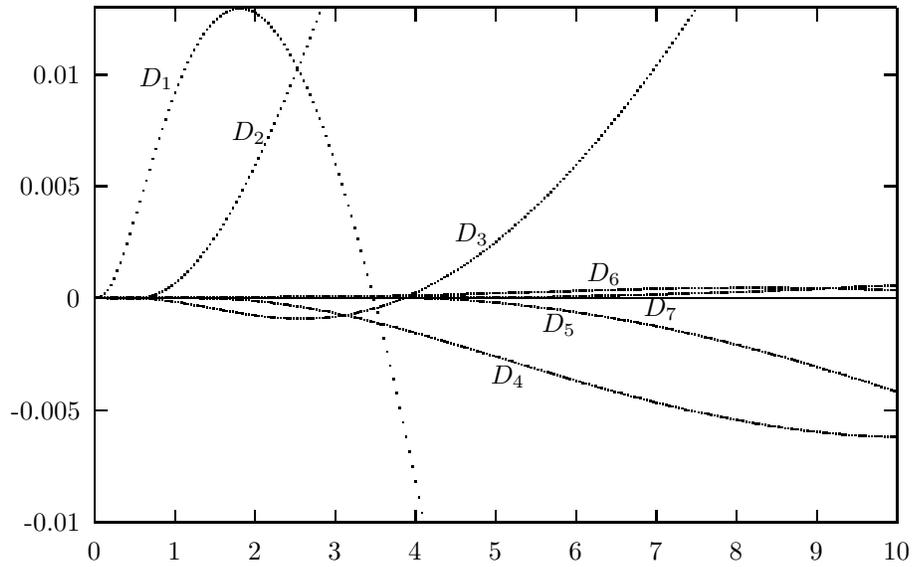

Figure 2: The differences $D_r(x) := A(x,r) - A(x)$, $r = 1, \ldots, 7$, for $x$ ranging in $[0, 10]$.

By letting $t = 1$, $s = 1/r$, and $\lambda = \sqrt{r}$ in Example 3.1, we get its analogue for (2.5) (note that $\phi(1/r, 1; \sqrt{r}) = \phi(1/\sqrt{r}, \sqrt{r}; 1)$). Everything else is unchanged from Example 3.1.

We also note that Monte Carlo approximation of $\phi(s, t; \lambda)$—which is very easy to implement—could have been useful if we did not have the incomparably more powerful approximation (2.16) at our disposal.

In order to demonstrate how extremely accurate the approximation (2.16) is, we need a powerful computing tool. We take this opportunity to recommend WIMS 'Function calculator', a module under the WWW Interactive Multipurpose Server (WIMS), available, for free, at

http://wims.unice.fr/wims/wims.cgi?module=tool/analysis/function.en.

Using WIMS function calculator, we calculated the difference between $A(x, r)$ and $A(x)$ (recall (3.1)) for various values of $x$ in $[0.01, 100]$ and $r$ in $[10, 100]$. The results are tabulated in Table 2. The $A(x)$ values were obtained effortlessly from the definition (2.7); however, in order to obtain the $A(x, r)$ values (from (2.8)), it was necessary to enter the function $\varphi(x, r)$ manually. From Table 2, it is evident that $A(x, r) - A(x)$ tends to 0 very rapidly as $r \to \infty$, uniformly for $x$ in bounded intervals $[0, a]$. This confirms our claim about the high quality of approximation (2.16), recalling, moreover, the last sentence in the first paragraph of this section.

*Remark* 3.1. It should be pointed out that the fact that $r$ is integer in Table 2 has no intrinsic significance. For example, the following rounded values of $A(x, r) - A(x)$ were obtained (cf. Table 2): $+5 \cdot 10^{-34}$ for $x = 0.01$, $r = 10.3$; $-8 \cdot 10^{-48}$ for $x = 0.10$, $r = 24.8$; $-1 \cdot 10^{-47}$ for $x = 1.00$, $r = 50.2$; and $-7 \cdot 10^{-46}$ for $x = 10.0$, $r = 99.7$.

Recall from Example 3.1 that we have implemented the computation of $\phi(s, t; \lambda)$ (Eq. (1.3)) in the Java language. In this setting, a significant drawback of Eq. (1.3) (in addition to that indicated in Example 3.1) is that the term $a_k := \lambda^k [t - (k-1)s]^k / k!$, assuming it is a double-precision number, would return the value Infinity if it overflows the value of `Double.MAX_VALUE`, which is approximately $1.8 \cdot 10^{308}$. (We note that $a_k$ should be computed as $\prod_{i=1}^{k} \lambda[t - (k-1)s]/i$, in order to avoid potential Infinity in the numerator or denominator.) As an example, if we set $s = 0.1$, $t = 500$, and $\lambda = 2$, then $a_k$ yields Infinity for $376 \leq k \leq 1431$; hence, $\phi(s, t; \lambda)$ would not be computed. On the other hand, $\phi^\star(s,t;\lambda)$ has been calculated, according to (2.16), to a very high precision (validated by comparison with WIMS). Moreover, changing $t$ in this example from 500 to 400 leads to the incorrect value 0 for $\phi(s,t;\lambda)$ (using Java), whereas, by comparison with the value obtained for $\phi^\star(s,t;\lambda)$ (using Java or WIMS), the correct value should be approximately $1 \cdot 10^{-54}$.

Table 3: Tabulated rounded values of $S_N(s,t;\lambda) - \bar{\phi}^\star(s,t;\lambda)$, $N = 3,\ldots,7$, for $s = 10^{-5}$, $t = 1$, and $\lambda = 100, 200$, and $400$.

| $\lambda$ | N | | | | |
|---|---|---|---|---|---|
| | 3 | 4 | 5 | 6 | 7 |
| 100 | $+5 \cdot 10^{-6}$ | $-1 \cdot 10^{-7}$ | $+2 \cdot 10^{-9}$ | $-4 \cdot 10^{-11}$ | $+5 \cdot 10^{-13}$ |
| 200 | $+1 \cdot 10^{-3}$ | $-9 \cdot 10^{-5}$ | $+7 \cdot 10^{-6}$ | $-4 \cdot 10^{-7}$ | $+2 \cdot 10^{-8}$ |
| 400 | $+2 \cdot 10^{-1}$ | $-7 \cdot 10^{-2}$ | $+2 \cdot 10^{-2}$ | $-5 \cdot 10^{-3}$ | $+1 \cdot 10^{-3}$ |

At this point, it is interesting to compare the approximation (2.16) with its counterpart (2.1) for the case when $t$ is fixed and $s$ is small. For this purpose and for later use, we shall consider the corresponding complements. We define

$$S_N(s,t;\lambda) = \sum_{n=1}^{N} (-1)^{n+1} P_n(\lambda t) \frac{\lambda^n s^n}{n!}, \qquad (3.2)$$

viewing it as a function of $s$. By (2.1),

$$\bar{\phi}(s,t;\lambda) = S_N(s,t;\lambda) + O(s^{N+1}) \quad \text{as } s \downarrow 0. \qquad (3.3)$$

We also define $\bar{\phi}^\star(s,t;\lambda) = 1 - \phi^\star(s,t;\lambda)$. We first recall from Example 3.1 and the paragraph that follows it that the exact formula (1.3) might be impractical when $r = t/s$ is very large, due to very long processing time. Thus (3.3) has an obvious merit over (1.3) in the asymptotic setting as $s \downarrow 0$, since $S_N$, like $\phi^\star$, is computed in 'no time'. However, as one would expect from (3.2) and Table 1, $N$ should be chosen with care when $\lambda$ is relatively large (especially if so is $t$), depending on the accuracy we wish to achieve. Hence the general superiority of $\bar{\phi}^\star$ over $S_N$ in approximating $\bar{\phi}$ (as $s \downarrow 0$). This point is demonstrated in Table 3 and Figure 3.

Having considered the case where $t/s$ is large—the important case in our context—we now claim that the approximation (2.16) may be effective already for $t$ larger than a few multiples of $s$, and, quite surprisingly, even in the entire range of interest $\{t \geq s\}$. This is apparent from Figure 2, combined with (3.1), for the case where $(r =) t/s = 1,\ldots,7$. As in Remark 3.1, we note that the fact that $r$ is integer in Figure 2 has no intrinsic significance. This point is demonstrated in Figure 4, showing that in the case $s = 1$, $\lambda = 1$ a quite satisfactory approximation is obtained uniformly for $t \geq 1$. Much more impressive approximations can be obtained by decreasing or (especially) increasing $\lambda$, as well as by decreasing or increasing $s$. We thus reveal another important advantage of approximation (2.16) over the exact form (1.3), namely, the simple structure of the former when $\phi^\star(s,t;\lambda)$ is viewed as a function of $t\,(\geq s)$. For example, by (1.1), the distribution function of the waiting time $\tau(s;\lambda)$ until two events of $N$ occur within $s$ time units can be approximated, for $t \geq s$, as

$$\mathrm{P}(\tau(s;\lambda) \leq t) \approx \bar{\phi}^\star(s,t;\lambda) = 1 - A(\lambda s)\mathrm{e}^{-\lambda[1-B(\lambda s)]t}.$$

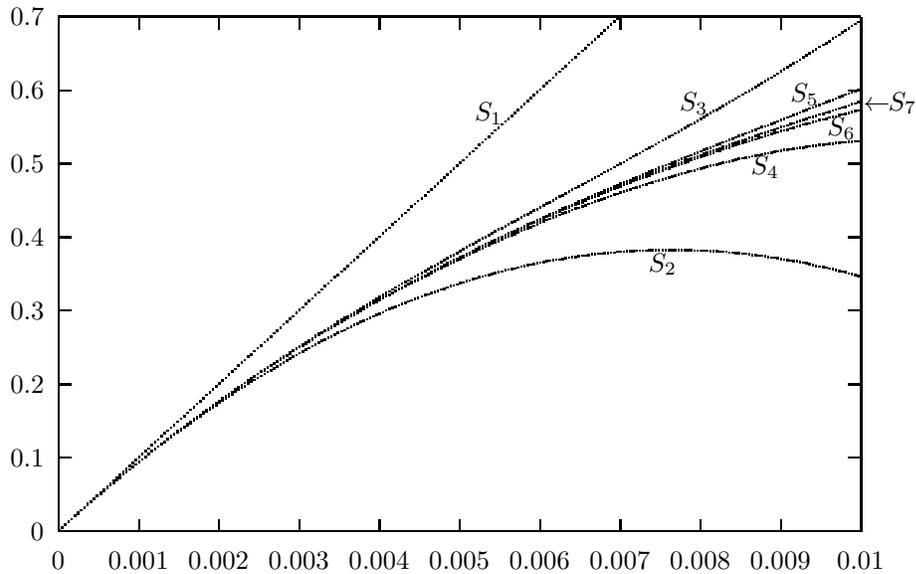

Figure 3: The first seven polynomials $S_N(s,t;\lambda)$ for $t = 1$, $\lambda = 10$, and $s$ ranging in $[0, 0.01]$.

(Note the resemblance to an exponential distribution function.)

# Appendix A:
# Comparison with Gates and Westcott's result

In this appendix we elaborate on the point made in the second paragraph of Section 1 regarding the fundamental result of this paper (namely, Theorem 2.2).

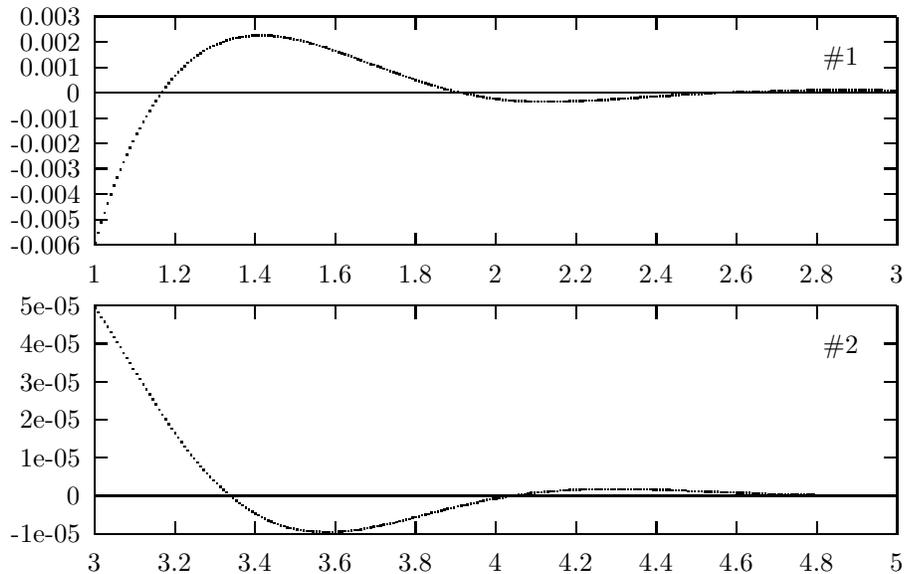

Figure 4: The difference $\phi^\star(s,t;\lambda) - \phi(s,t;\lambda)$ for $s=1$, $\lambda=1$, and $t$ ranging in $[1,3]$ (#1) and $[3,5]$ (#2).

In our notation, Gates and Westcott [3], who considered the distribution of the scan statistic $Y_N(s,t)$, gave the following expression for the function $\varphi(s,t;\lambda) = e^{\lambda t}\phi(s,t;\lambda)$ [3, Eq. (2.3)]:

$$\varphi(s,t;\lambda) = \frac{\exp[p(t+s)]}{1+sp} + \exp(pt)o(1) \quad \text{as } t \to \infty, \qquad (A.1)$$

where $p = p(\lambda, s)$ is the solution of $p\exp(sp) = \lambda$. It is readily checked that the solutions $p(\lambda, s)$ and $B(\lambda s)$ are related by $p(\lambda, s) = \lambda B(\lambda s)$ and, in turn, that (A.1) is equivalent to (2.11). However, our presentation is evidently more elegant and instructive.

Much more importantly, Gates and Westcott underestimated the accuracy and validity of the approximation. They considered it "extremely accurate (typically to six or seven significant figures) even for small $L$", where $L$ is $t$ in our notation, and illustrated its accuracy in Table 1 of their paper. However, the values tabulated in that table are relatively inaccurate and misleading, the most notable example being the values 0.5439808 and 0.5439807 for $\phi(1, 20; 0.2)$ and $\phi^\star(1, 20; 0.2)$, respectively, whereas using WIMS we have found that $\phi(1, 20; 0.2) - \phi^\star(1, 20; 0.2) \approx 9.7 \cdot 10^{-33}$. As for the validity of the approximation, Gates and Westcott apparently did not realize that the approximation is valid for large values of $t/s$ in general, and in particular in the asymptotic setting as $s \downarrow 0$; their Table 1 presents, in addition to the approximation of $\phi(1, 20; 0.2)$, approximations to $\phi(1, 20; 0.3)$ and $\phi(1, 10; 1)$. That table also presents the corresponding values of Naus's approximation mentioned in Section 1; roughly, the errors from the exact $\phi$ values are of order $10^{-4}$.

A final note is in order concerning the derivation of (A.1) by Gates and Westcott. In our notation, Gates and Westcott noted from Penrose and Elvey [6] that the Laplace transform of $\varphi(s, t; \lambda)$, as a function of $t$, is given by

$$\int_0^\infty e^{-zt} \varphi(s, t; \lambda) \, dt = \frac{1}{z} + \frac{\lambda}{z^2[1 - \lambda \exp(-sz)/z]} \tag{A.2}$$

(say, for all $z$ satisfying $\operatorname{Re}(z) > \lambda$), and that inverting it by residues gives (A.1). However, justifying the second part of this statement rigorously might be a little inconvenient, especially for one who is not familiar with the calculus of residues. The delicate point is that the right-hand side of (A.2), when extended to a meromorphic function on the whole complex plane, has infinitely many poles. (In this context, see sections 3 and 4 in Penrose and Elvey, where a somewhat more general case is treated.) In this respect, application of the Wiener–Ikehara theorem is a good alternative.

In any case, Gates and Westcott's result has not received its due recognition.

# Appendix B:
# The Wiener–Ikehara theorem

For the convenience of the reader, we quote Theorem 4.2 from Korevaar [4].

**Theorem B.1** (The Wiener–Ikehara theorem). *Let $S(t)$ vanish for $t < 0$, be nonnegative and nondecreasing for $t \geq 0$, and such that the Laplace–Stieltjes transform*

$$f(z) = \mathcal{L}dS(z) = \int_{0^-}^\infty e^{-zt} \, dS(t) = z \int_0^\infty S(t) e^{-zt} \, dt, \quad z = x + iy,$$

*exists for $\operatorname{Re}(z) = x > 1$. Suppose that for some constant $A$, the analytic function*

$$g(z) = f(z) - \frac{A}{z-1}, \quad x > 1,$$

*has a boundary function $g(1 + iy)$ in the following sense. For $x \downarrow 1$, the function $g_x(iy) = g(x + iy)$ converges to $g(1 + iy)$ uniformly or in $L^1$ on every finite interval $-\lambda < y < \lambda$. Then,*

$$e^{-t} S(t) \to A \quad \text{as } t \to \infty.$$

## Acknowledgements

I am still indebted to Mr and Mrs Shapack for their generous support during my PhD studies (2005–2008).

# Erratum/Addendum to "On probabilities associated with the minimum distance between events of a Poisson process in a finite interval"

Shai Covo

Department of Mathematics

Bar Ilan University

52900 Ramat-Gan, Israel

E-mail: green355@netvision.net.il

November 1, 2018

## Abstract

In this addendum, we further discuss numerical issues concerning the computation of the probability denoted $\phi(s,t;\lambda)$ in [1]. In particular, and most importantly, we correct the naive claim made in the abstract of [1] that the explicit expression of $\phi(s,t;\lambda)$ becomes intractable for very large values of $t/s$; rather, we show that it may not be applicable for large values of $\lambda t$.

We use the notation of [1]. At first glance, the expression on the right of

$$\phi(s,t;\lambda) = e^{-\lambda t}\left(1 + \sum_{k=1}^{\lceil t/s \rceil} \frac{\lambda^k [t-(k-1)s]^k}{k!}\right) \quad \text{(A.1)}$$

[1, Eq. (1.3)] seems intractable for very large values of $t/s$. In some respects this impression is correct (comparing to the approximations given in [1]). However, in order to reveal a significant practical drawback of (A.1), we should rather turn our attention to large values of $\lambda t$. To exemplify this claim, we suppose—as in [1]—that the computation of $\phi(s,t;\lambda)$ is implemented in Java (using double precision). In this setting, the value of an expression is Infinity if it overflows the value of `Double.MAX_VALUE`, which is approximately $1.8 \cdot 10^{308}$. Thus, the term $a_k := \lambda^k[t-(k-1)s]^k/k!$ should be computed with some care, say as $a_k = \prod_{i=1}^{k} \lambda[t-(k-1)s]/i$. Doing so, no overflow problem would arise as long as $\lambda t$ (rather than $t/s$) is not too large. Indeed, since $\phi(s, t; \lambda)$ is a probability, the expression multiplying $\mathrm{e}^{-\lambda t}$ on the right of (A.1) is bounded from above by $\mathrm{e}^{\lambda t}$; so are the partial products of $a_k$, since $\prod_{i=1}^{j} \lambda[t - (k-1)s]/i \leq (\lambda t)^j/j!$ ($1 \leq j \leq k$). Thus, the condition $\lambda t < c^*$ where $c^* = \ln(\texttt{Double.MAX\_VALUE}) \approx 709.78$ guarantees the applicability of (A.1). (It should be noted here that $\texttt{Double.MIN\_VALUE} = 4.9 \times 10^{-324} < \mathrm{e}^{-744}$.) For larger values of $\lambda t$, the story is quite different. Suppose that $\lambda t > c^*$. Then, in view of (A.1), an overflow problem will occur if $\phi(s, t; \lambda) > p^*$ where $p^* = \mathrm{e}^{-\lambda t}\texttt{Double.MAX\_VALUE}$. (For example, if $\lambda t = 716$, then $p^* \approx 2 \cdot 10^{-3}$.) In order to mitigate this problem, we first rewrite (A.1) as follows:

$$\phi(s, t; \lambda) = \mathrm{e}^{-\lambda t} + \sum_{k=1}^{\lceil t/s \rceil} \prod_{i=1}^{k} \exp(-\lambda t/k) \frac{\lambda[t - (k-1)s]}{i}. \tag{A.2}$$

As one would expect, this representation allows for substantially larger values of $\lambda t$ than previously. However, the improvement is not drastic, as we now exemplify. Let $r = 4 \times 10^6$, and set $s = 1/\sqrt{r}$, $t = \sqrt{r}$, and $\lambda = 1$. Hence, $\lambda t = 2000$. The approximation [1, Eq. (2.16)]

$$\phi^\star(s, t; \lambda) := A(\lambda s)\mathrm{e}^{-\lambda[1 - B(\lambda s)]t} \tag{A.3}$$

of $\phi(s, t; \lambda)$ gives us $\phi(s, t; \lambda) = 0.36815525...$ (in agreement with the estimate $\mathrm{e}^{-1} + 1.5\mathrm{e}^{-1}/\sqrt{r} = 0.36815535...$ indicated in [1, Example 3.1]). Truncating the sum in (A.2) at $k = 1942$ contributes less than 11% to the value of $\phi(s, t; \lambda)$; the 1943th term, on the other hand, contributes Infinity. (Hence (A.2) is useless here.) To understand this phenomenon, consider the partial products of the $k$th term:

$$\Pi(j; k) := \prod_{i=1}^{j} \exp(-\lambda t/k) \frac{\lambda[t - (k-1)s]}{i}, \quad 1 \leq j \leq k.$$

While, theoretically, $\Pi(k; k)$ is eventually small, it might reach Infinity on the way. In our example, the theoretical value of $\Pi(1943; 1943)$ is approximately $0.00155$ (found by (A.4), below) but $\Pi(699; 1943) = \text{Infinity}$; this results in $\Pi(1943; 1943)$ (and hence also $\phi(s, t; \lambda)$) being equal to Infinity. See also Figure 1. It is thus natural to consider the following modification of (A.2):

$$\phi(s, t; \lambda) = \mathrm{e}^{-\lambda t} + \sum_{k=1}^{\lceil t/s \rceil} \prod_{i=1}^{k} \exp(-\lambda t/k) \frac{\lambda[t - (k-1)s]}{(k - i + 1)}. \tag{A.4}$$

However, it turns out that computing $\phi(s, t; \lambda)$ this way does not lead to an essential improvement, if any.

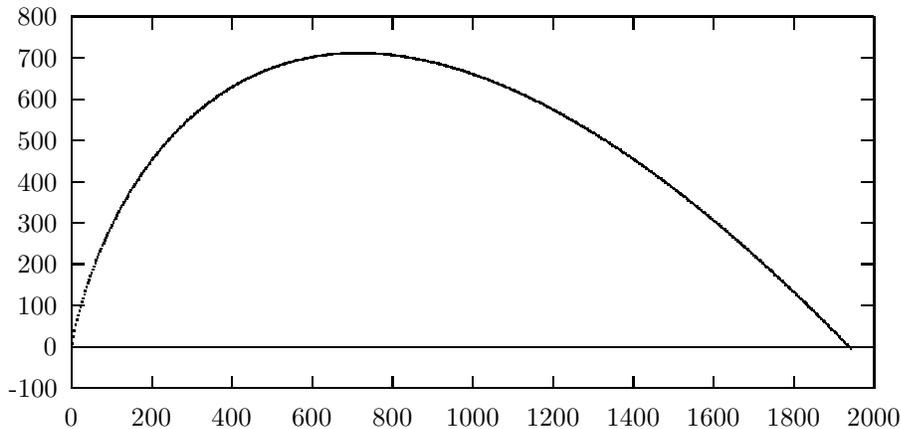

Figure 1: Plot of the function $\ln \Pi(j;k)$, $k = 1942$, $j = 1,\ldots,k$. Here $\max_j \ln \Pi(j;k) \approx 709.56$, just below the permitted maximum of $c^* \approx 709.78$.

Having partially resolved the computational overflow problem of (A.1), we now turn to consider the computational time problem indicated in [1, Example 3.1]. For the purposes of discussion, we assume that $\phi(s,t;\lambda)$ is computed using (A.1) and, of course, that $t/s$ is very large; in view of the preceding discussion, we assume that $\lambda t$ is not too large. In particular, we assume that $\lambda t$ is substantially smaller than $t/s$; equivalently, $\lambda s$ is substantially smaller than 1. (In view of (A.6), this restriction is very mild.) In this setting, the aforementioned problem is resolved by the following observation. For any $a > 0$ fixed and $\lambda t$ sufficiently large (depending on $a$),

$$(0 <) \phi(s,t;\lambda) - e^{-\lambda t} \sum_{|k-\lambda t| < a\sqrt{\lambda t}} \frac{\lambda^k [t-(k-1)s]^k}{k!}$$
$$< \sum_{|k-\lambda t| \geq a\sqrt{\lambda t},\ k \geq 0} e^{-\lambda t} \frac{(\lambda t)^k}{k!} \approx 2(1 - \Phi(a)), \quad \text{(A.5)}$$

where $\Phi$ denotes the distribution function of the $N(0,1)$ distribution. The approximation on the right follows by a straightforward application of the central limit theorem. (The basic idea is that if $X_i$ are i.i.d. Poisson(1) and $\lambda t$ is an integer, then $\sum_{i=1}^{\lambda t} X_i$ is Poisson($\lambda t$).) It is very rough, as one can conclude by carrying out the proof, but still provides a satisfactory order of magnitude. The implication of (A.5) is that $\phi(s,t;\lambda)$ can be approximated very well by summing a relatively small number of terms that dominate the calculation. As an example, let $s = 1/400$, $t = 400$, and $\lambda = 1$. (Hence, $\phi(s,t;\lambda) \approx 0.369$.)

Approximating $\phi(s,t;\lambda)$ according to (A.5) with $a = 5$ (hence, 199 terms are taken) gives an error of approximately $2.24 \times 10^{-7}$, in agreement with the hypothetical upper bound $2(1 - \Phi(5)) \approx 5.73 \times 10^{-7}$. Increasing $a$ to 6 (hence, 239 terms are taken) gives an error of approximately $1.07 \times 10^{-9}$, in agreement with the hypothetical upper bound $2(1 - \Phi(6)) \approx 1.97 \times 10^{-9}$. As a further example, let $s = 10^{-3}$, $t = 10^3$, and $\lambda = 0.4$. (Hence, $\phi(s,t;\lambda) \approx 0.852$.) Approximating $\phi(s,t;\lambda)$ by summing from $k = 190$ to $650$ in (A.1), gives an error of approximately $5.77 \times 10^{-31}$ (confirmed, using WIMS "Function calculator", by comparison with the value of $\phi^\star(s,t;\lambda)$). We conclude that the apparent computational time problem indicated in [1, Example 3.1] is simply resolved by truncating the sum in (A.1) according to (A.5).

We take this opportunity to make the following remark.

*Remark* A.1. Extending Remark 1.1 of [1], we point out that $\phi$ depends on $s$, $t$, and $\lambda$ only through $\lambda s$ and $\lambda t = \lambda s (t/s)$. In particular, we have that

$$\phi(s,t;\lambda) = \phi(\lambda s, \lambda t; 1) = \phi(1, t/s; \lambda s). \tag{A.6}$$

The new representation, $\phi(s,t;\lambda) = \phi(\lambda s, \lambda t; 1)$, is instructive in view of (A.3), (A.5), (A.7), and [1, Eqs. (2.4) and (2.5)]. Moreover, it shows that the computational problem indicated in the first paragraph cannot be simply resolved by some transformation of the parameters. On the other hand, the representation $\phi(s,t;\lambda) = \phi(1, t/s; \lambda s)$ arises from the setting of Theorem 2.2 of [1]—the fundamental result of the paper.

Suppose that (A.1) (including its variants, (A.2) and (A.4)) is not applicable due to a large value of $\lambda t$. If $\lambda s \cdot \lambda t$ is smaller than 1 or so (the smaller the better), then $\phi(s,t;\lambda)$ can be approximated fairly accurately by the asymptotic result [1, Eq. (2.1)]

$$\phi(s,t;\lambda) - \sum_{n=0}^{N} (-1)^n P_n(\lambda t) \frac{\lambda^n s^n}{n!} = O(s^{N+1}) \quad \text{as } s \downarrow 0. \tag{A.7}$$

(The larger $N$ is, the better the approximation.) However, such approximation is inferior to the outstanding one (A.3). This point is illustrated in Table 1 (cf. also [1, Table 3]). For convenience of implementation, it pays to rearrange the expression for $P_n(x)$ [1, Section 2] as follows:

$$P_n(x) = \sum_{j=0}^{n} \left[ \sum_{i=j}^{n} \binom{n}{i} (n-1)^{n-i} S(i,j) \right] x^j$$

(where $S(i,j)$ denote the Stirling numbers of the second kind). Finally, it should be stressed that (A.7) is not effective when $\lambda s \cdot \lambda t$ is substantially greater than 1, say greater than 3. In this case, and in general when $\lambda t$ is large, a satisfactory approximation to $\phi(s,t;\lambda)$ may be obtained simply by the asymptotic estimate $e^{-\lambda^2 st}$ (cf. [1, Eqs. (2.4) and (2.5)]). For example, if $\lambda s = 10^{-3}$ and $\lambda t = 4 \times 10^3$, then $\phi(s,t;\lambda) - e^{-\lambda^2 st} \approx 1.1 \times 10^{-4}$; if $\lambda s = 10^{-5}$ and $\lambda t = 4 \times 10^5$, then $\phi(s,t;\lambda) - e^{-\lambda^2 st} \approx 1.1 \times 10^{-6}$.

Table 1: Tabulated rounded values of $\phi^\star(s,t;\lambda) - \sum_{n=0}^{N} (-1)^n P_n(\lambda t)\lambda^n s^n/n!$, $N = 6,\ldots,10$, for $\lambda = 1$, $t = 10^4$, and $s = 10^{-4}$, $10^{-6}$, and $10^{-8}$.

| $s$ | $N$ | | | | |
|---|---|---|---|---|---|
| | 6 | 7 | 8 | 9 | 10 |
| $10^{-4}$ | $-1.8 \cdot 10^{-4}$ | $+2.2 \cdot 10^{-5}$ | $-2.5 \cdot 10^{-6}$ | $+2.6 \cdot 10^{-7}$ | $-2.3 \cdot 10^{-8}$ |
| $10^{-6}$ | $-2.0 \cdot 10^{-18}$ | $+2.5 \cdot 10^{-21}$ | $-2.8 \cdot 10^{-24}$ | $+2.8 \cdot 10^{-27}$ | $-2.5 \cdot 10^{-30}$ |
| $10^{-8}$ | $-2.0 \cdot 10^{-32}$ | $+2.5 \cdot 10^{-37}$ | $-2.8 \cdot 10^{-42}$ | $+2.8 \cdot 10^{-47}$ | $-2.5 \cdot 10^{-52}$ |